\documentclass[10pt,twocolumn,letterpaper]{article}

\usepackage{cvpr}
\usepackage{times}
\usepackage{epsfig}
\usepackage{graphicx}
\usepackage{amsmath}
\usepackage{amssymb}
\usepackage{amsfonts}
\usepackage{textcomp}
\usepackage{amsthm}

\usepackage[breaklinks=true,bookmarks=false]{hyperref}

\cvprfinalcopy 

\def\cvprPaperID{****} 

\newcommand{\qq}{{\bf x}}

\begin{document}

\title{Signal Classification in Quotient Spaces via \\ Globally Optimal Variational Calculus}

\author{Gregory S. Chirikjian\\
Laboratory for Computational Sensing and Robotics \\
Johns Hopkins University\\
{\tt\small gchirik1@jhu.edu}
}

\maketitle

\begin{abstract}
A ubiquitous problem in pattern recognition is that of matching an observed time-evolving pattern (or signal) to a gold standard
 in order to recognize or characterize the meaning of a dynamic phenomenon. Examples include matching sequences of images in two videos, matching audio signals in speech recognition, or matching framed trajectories
in robot action recognition. This paper shows that all of these problems can be aided by
reparameterizing the temporal dependence of each signal individually to a universal standard timescale that
allows pointwise comparison at each instance of time. Given two sequences, each with $N$ timesteps,
the complexity of the algorithm has a cost of $O(N)$, which is an improvement on the most common method for matching two signals,
i.e., dynamic time warping. The core of the approach presented here is that the universal standard timescale is found
by solving a variational calculus problem in which the cost functional reflects the amount of change that takes
place as measured in the original temporal variable, and then produces a mapping to a new temporal variable in
which the amount of change is globally minimized. The result builds on known facts in differential geometry.
\end{abstract}

\section{Introduction}

Consider two time-evolving sequences, or signals, $X_1(t)$ and $X_2(t)$ which can be scalar, vector, matrix, or Lie group quantities of the same type.
Without loss of generality, let $t \in [0,1]$. Let $S$ denote the space in which all such signals evolve.
Then $X(t)$ can be viewed as a map
$$ X:[0,1] \,\rightarrow S $$
and the space of all signals is $[0,1] \times S$.

Suppose that there exists a metric (distance function)
$d:S \times S \,\rightarrow\, \mathbb{R}_{\geq 0}$, thereby making $(S,d)$ a metric space. In general it can be the case that $d(X_1(t), X_2(t))$ will not be small even if $X_1(t)$ and $X_2(t)$ represent
fundamentally the same dynamic phenomenon. This can happen for two reasons: (1) The sequences can have a different temporal evolution along the fundamentally same path; (2) Nuisance parameters such as perspective, background noise, or
signal decimation can cloud the underlying similarity. The first of these problems can be addressed by considering the
internal (temporal) dependence of the signal that act on the time interval $[0,1]$, and the second can be described by external transformations that act on the space $S$. In this paper both of these are considered, as well as joint transformations that act on the whole space $[0,1] \times S$ in a coupled way.

 As a first example, consider when $X_i(t)$ are two scalar functions each describing the audio signal of spoken text ``The rain in Spain stays mainly in the plain.''\footnote{See
YouTube for this part of the movie {\it My Fair Lady}} $X_1(t)$ could be the template of how this phrase should be spoken, and $X_2(t)$ could be how someone with an accent says (or sings) the same phrase. Though $d(X_1(t), X_2(t))$ will not
be small if the second person has a strong accent or carries certain syllables longer while singing, the expansion or
contraction of certain syllables over time can be compensated for by instead reparameterizing both to a standard timescale by defining
\begin{equation}
Y_i(t) \doteq X_i(\tau_i(t))
\label{Ydef}
\end{equation}
where
$$ \tau_i: [0,1] \,\rightarrow\, [0,1] $$
are smooth monotonically increasing functions with smooth inverse. The set of all such functions forms a group $({\cal T}, \circ)$ under the operation of
composition of functions. That is, given $\tau_1, \tau_2 \in {\cal T}$,
$$ (\tau \circ \tau')(t) \,\doteq\, \tau(\tau'(t)) $$
is also in ${\cal T}$, and satisfies all of the group axioms such as associativity, and the inverse group element $\tau^{-1}(t)$ is the inverse of the function $\tau(t)$
which exists due to monotonicity, and the identity element is $e(t) = t$. Let us call this the {\it temporal reparameterization group, TRG}. This is
actually an infinite-dimensional group of transformations that act on $[0,1]$.
Then we can say that $X_1(t)$ and $X_2(t)$ are fundamentally the same if there exist $\tau_1, \tau_2 \in {\cal T}$ such that\footnote{As a practical
matter, often one seeks to minimize the integral of the square of a metric to eliminate square roots under the integral.}
\begin{equation}
\int_{0}^{1} [d(Y_1(t), Y_2(t))]^2 dt \,\approx\, 0
\label{}
\end{equation}
even if $d(X_1(t), X_2(t)) \gg 0$ for all values of $t$.

Of course, it would not be feasible to construct a search over the space ${\cal T} \times {\cal T}$ since ${\cal T}$ is infinite dimensional.
Herein lies one of the fundamental contribution of this paper: It is possible to independently obtain $\tau_i(t)$ resulting in reparameterizations $Y_i(t)$
each on a {\it universal standard timescale, UST}, using a particular variational calculus formulation that is realizable in $O(N)$ computations where
$N$ is the number of time steps in the recorded sequence.

This is not limited to the scalar case described in the audio example. For example, if a robot arm is doing free-form manufacturing with a milling tool as its end effector, it could be that the same path is traversed by two different trajectories implementing the same task, but with different dwell times or different rates along different parts of the trajectory. In this case the problem might be modelled with $S = SE(3)$, which is the 6-dimensional Lie group of
rigid-body displacements with group operation being matrix multiplication when elements are expressed as matrices of the form
\begin{equation}
X = \left(\begin{array}{ccc}
R && {\bf t} \\ \\
{\bf 0}^T && 1 \end{array} \right).
\label{se3def}
\end{equation}
Here $R$ is a $3\times 3$ rotation matrix, ${\bf t}$ is a 3D translation vector, and ${\bf 0}^T$ is a row of zeros.
Then a metric (distance function) on this space is
\begin{equation}
d(X_1, X_2) = \|\log(X_{1}^{-1} X_{2})\|_F
\label{ddddef}
\end{equation}
where $\|X\|_F = \sqrt{{\rm tr}(XX^T)}$ is the Frobenius matrix norm. Note that $d(X_1, X_2)$ is left invariant in the sense that
$d(X_0 X_1, X_0 X_2) = d(X_1, X_2)$ for arbitrary $X_0,X_1,X_2 \in SE(3)$.

Group theory appears in several ways in this problem.
In additional to $S$ possibly being a group in some contexts (such as when the signal is a robot trajectory), and the
TRG is a group which is used to internally quotient out the effects of temporal fluctuations, it can be desirable to simultaneously quotient out the effects of other ``nuisance groups'' \cite{Soatto} that externally act on $S$. For example, if two video sequences of a person waving is presented
from two different viewpoints, and we wish to discern whether the actions are actually those of waving as opposed to throwing a ball,
then the the effects of viewpoint on each image in each of the sequences can be quotiented out as well. For example, if $G$ is a group such as
$SL(3,\mathbb{R})$ when considering homographies \cite{Mahony}, or $SE(2)$ when considering rigid template matching in the image plane \cite{dover}, or $SO(3)$ when considering image matching from fish-eye lenses \cite{13Makadia2010}, or an affine transformation approximating a perspective transformation in a pinhole camera
\cite{Hartley}, then \cite{dover}
\begin{equation} D_G(X_1, X_2) \doteq \min_{g \in G} d(X_1, g \cdot X_2)
\label{Ddef}
\end{equation}
is a metric on the quotient space $G\backslash S$ where $\cdot$ denotes the action of $G$ on $S$.
Invariant recognition of signals then becomes one of matching in the double quotient space $G\backslash S/{\cal T}$.
The contribution of this paper is developing the mathematical framework to do this. But first, a brief review of
what is usually done in the literature is provided.

\subsection{Related Literature}

The pattern recognition literature is immense, and is divided into the subcommunities centered around different application areas
such as computer vision and image understanding, speech recognition, and robotics.
In all of these areas deep learning has made tremendous strides in recent years. See, for example, \cite{LeCun,Schmidhuber,Lee}.
A method used in information processing and pattern recognition, originally developed for audio signals
\cite{Rabiner,Berndt}, is {\it dynamic time warping}. In the simplest implementation of this method, a measure of similarity between
each point in two sequences is used to generate a pairing cost in a bipartite graph. The resulting matching is usually computed in $O(N^2)$ time when
there are $O(N)$ points in each sequence, though algorithms exist to reduce this computational burden somewhat \cite{Salvador}.

This paper provides the mathematical framework for a very different alternative. Rather than morphing (or warping) one audio or image
sequence to fit another, the goal here is to reparameterize each time-varing object to its own natural time scale, and to simultaneously
quotient out the effects of nuisance groups. (Much of this paper is concerned with the theoretical underpinnings justifying how to do these computation). Then two
sequences can be compared directly pointwise. The genesis of this idea (without the nuisance groups) was a discussion on this topic in the context of a particular application with Ms. Yixin Gao \cite{Yixin}.

\subsection{Structure of the Remainder of the Paper}

In Section \ref{sec1} a class of problems in the Calculus of Variations that is directly related to efficiently
selecting elements of the temporal reparameterization group, ${\cal T}$, for establishing correspondences between signals in general
Riemannian metric spaces $(S, d)$ is articulated. In particular, a proof is given that this
class of variational calculus problems yields a globally optimal solution. Section \ref{sec2} then discusses
more widely the question of when variational calculus problems are guaranteed to have globally optimal solutions
generated by the Euler-Lagrange equation, and how a bootstrapping procedure can be used to expand the scope of problems
that have such globally optimal solution.

\section{Global Optimality in Variational Calculus to Reduce Searches Over $S$ to $S/{\cal T}$} \label{sec1}

The Calculus of Variations addresses the problem of seeking vector-valued functions ${\bf x}(t)$ that extremize
functionals of the form
\begin{equation}
J = \int_{0}^{1} f\left({\bf x},\dot{\bf x},t\right)\,dt
\label{Jdefapp31245}
\end{equation}
where $ \dot{\bf x} = d{\bf x}/dt$.
As in usual calculus, the result can be either saddle-like solution or local or global minimum, maximum.
Necessary conditions for such solutions are given by the Euler-Lagrange equations:
\begin{equation}
\frac{\partial f}{\partial {\bf x}} -
\frac{d}{dt} \left(\frac{\partial f}{\partial \dot{\bf x}} \right) = {\bf 0}
\label{eulerlag}
\end{equation}
where derivatives with respect to vectors are interpreted as gradients. When ${\bf x}(t)$ is one dimensional it is denoted as $x(t)$.

The Euler-Lagrange equations only provide ``first order'' necessary conditions
for a local (or weak) extremum. Stronger necessary condition due to Jacobi also exist, but even then
there is usually no guarantee that a solution
of the Euler-Lagrange equations will
be globally optimal. However, in certain situations (including optimal temporal reparameterization), the structure of the function
$f(\cdot)$ will guarantee that the solution generated by the
Euler-Lagrange equations is in fact a globally optimal solution.
In particular, we have the following.
\\ \\
\noindent
{\bf THEOREM 1}:
{\it When the integrand in the cost functional (\ref{Jdefapp31245}) is of the form
\begin{equation}
f(x,\dot{x}) = \dot{x}^2 \mathfrak{g}(x)
\label{optcost}
\end{equation}
where $\mathfrak{g}(x)$ is differentiable and $\mathfrak{g}(x(t)) > 0$ for all values of $t \in [0,1]$,
then the solution generated by (\ref{eulerlag}) subject to the boundary conditions $x(0) = 0$ and $x(1) =1$ is globally minimal}.
\footnote{Such cost functions arise in reparameterization in natural ways.}

\begin{proof}
Evaluating (\ref{eulerlag}) with (\ref{optcost}) gives
\begin{equation}
 2 \ddot{x} \mathfrak{g} + \dot{x}^2 \frac{\partial \mathfrak{g}}{\partial x} = 0.
 \label{origode}
\end{equation}
Multiplying both sides by $\dot{x}$ and integrating yields the
exact differential
$$ \frac{d}{dt}(\dot{x}^2 \mathfrak{g}) = 0. $$
Integrating both sides with respect to $t$ and isolating $\dot{x}$ yields
$$ \dot{x} = c\, \mathfrak{g}^{-\frac{1}{2}}(x) $$
where $c$ is the arbitrary constant of integration. With the
boundary conditions $x(0) =0$ and $x(1) =1$, we can then
write
$$ F(x^*) \doteq \frac{1}{c} \int_{0}^{x^*} \mathfrak{g}^{\frac{1}{2}}(\sigma)\,d\sigma = t, $$
where
$$ c = \int_{0}^{1} \mathfrak{g}^{\frac{1}{2}}(\sigma)\,d\sigma. $$
The notation $x^*$ indicates that this is the unique solution obtained from the
Euler-Lagrange equations that satisfies the boundary conditions.

The function  $F(x^*) = t$ can be
inverted ($F$ is monotonically increasing since
$\mathfrak{g}(x) > 0$) to yield $x^* = F^{-1}(t).$

To see that this solution is globally optimal, substitute
\begin{equation}
\dot{x}^* = \mathfrak{g}^{-\frac{1}{2}}(x^*) \int_{0}^{1} \mathfrak{g}^{\frac{1}{2}}(\sigma)\,d\sigma
\label{soln}
\end{equation}
into the cost functional
$$ J(y) = \int_{0}^{1} \mathfrak{g}(y) \dot{y}^2 dx $$
where $y(t)$ is \underline{any} function in ${\cal T}$.
Then
$$ J(x^*) = \left(\int_{0}^{1} \mathfrak{g}^{\frac{1}{2}}(x^*)\,dx^* \right)^2
= \left(\int_{0}^{1} \mathfrak{g}^{\frac{1}{2}}(y)\,dy \right)^2, $$
where the second equality is simply a change of name of the dummy variable of integration.
Furthermore, since $x^*$ and $y$ are both functions of time, we can change the domain of integration as
$$  J(x^*) = \left(\int_{0}^{1} \mathfrak{g}^{\frac{1}{2}}(y(t)) \,\dot{y}\, dt \right)^2. $$
Since in general, from the Cauchy-Schwarz inequality,
$$ \left(\int_{0}^{1} f(t)\,dt \right)^2 \,\leq \int_{0}^{1} [f(t)]^2 dt, $$
we see that by letting $f(t) = \mathfrak{g}^{\frac{1}{2}}(y) \dot{y}$ that
$$ \left(\int_{0}^{1} \mathfrak{g}^{\frac{1}{2}}(y(t)) \,\dot{y}\, dt \right)^2 \,\leq\,
\int_{0}^{1} \mathfrak{g}(y) (\dot{y})^2 dt $$
and hence
$$ J(x^*) \leq J(y) $$
where $x^*(t)$ is the solution generated by the Euler-Lagrange
equation and $y(t)$ is any function in ${\cal T}$. Therefore $x^*(t)$ is a globally minimal solution.
\end{proof}
Note: The condition that $\mathfrak{g}(x)$ is differentiable was required to use the Euler-Lagrange equation, but
if (\ref{soln}) had been hypothesized independently, this condition could be relaxed to continuity (which
is still required in order to be able to invert $F(\cdot)$), and
the global optimality of the solution would persist.

As an example of how (\ref{optcost}) arises, consider when $S = \mathbb{R}^{n\times m}$ and the metric $d$ is a matrix norm of the difference of two elements. Then
$$ d(X(t+dt), X(t)) = \|X(t+dt) - X(t)\| = \|dX/dt\|\, dt. $$
Recall that $Y(t) = X(\tau(t))$ is the reparameterized version of $X(t)$.
Then minimizing the integral over $t \in [0,1]$ of
$$ \left\|\frac{dY}{dt}\right\|^2 \,=\, \left\|\frac{d}{dt} X(\tau(t)) \right\|^2 = \left\|X'(\tau) \dot{\tau} \right\|^2 =
\left\|X'(\tau) \right\|^2 \dot{\tau}^2 $$
gives an $f(\cdot)$ of the form in (\ref{optcost}), with $\tau$ taking the place of $x$ and $\mathfrak{g}(\tau) = \left\|X'(\tau) \right\|^2$. (Here, of course, $X'(\tau) = dX/d\tau$.)
In other words, identifying the element of the TRG that {\it optimally} reparameterizes time
boils down to precisely solving the globally optimal variational calculus problem addressed in the above theorem.
The global optimality is critical because it means that there is a unique way to reparameterize the temporal dependence of a signal
so that the temporal fluctuations are spread out as evenly as possible.

\section{Bootstrapping Global Optimality to Larger Spaces} \label{sec2}

The theorem presented in the previous section begs the more general question of when solutions to variational problems are globally optimal.
To the author's knowledge two classes of such problems have been addressed in the literature.

First, in Riemannian geometry, a cost function of the form
$$ f({\bf x}, \dot{\bf x}) = \sqrt{\dot{\bf x}^T G({\bf x}) \dot{\bf x}} $$
where $T$ denotes the transpose of a vector or matrix and
$G({\bf x})$ is the metric tensor for a Riemannian manifold with negative sectional curvatures is know to have unique
geodesics, which hence globally minimize the functional (\ref{Jdefapp31245}). For example, in the Poincar\'{e} solid
n-dimensional open unit ball model of the hyperbolic space, the
metric tensor $G({\bf x}) = [g_{ij}({\bf x})]$ with ${\bf x} \in \mathbb{B}^n \subset \mathbb{R}^n$ is
$$ g_{ij}({\bf x}) = \frac{4 \delta_{ij}}{\left(1 - {\bf x}^T{\bf x}\right)^2} $$
where $\delta_{ij}$ is the Kronecker delta. This $G({\bf x})$
is known to have constant negative sectional curvature of value $-1$, thus guaranteeing that any geodesic connecting
two points has minimal length \cite{Ratcliffe,Cheeger,Benedetti,Eberlein}. In contrast, for a space of non-negative curvature such as the torus or sphere, geodesics exist that are not necessarily minimal length (e.g., one can take the long way around a great arc
to connect two points), and hence global minimality of length is not guaranteed for geodesics between arbitrary points in general Riemannian manifolds.

The second class of functions that globally minimize
(\ref{Jdefapp31245}) have been studied in the variational calculus literature. These are cost functions $f({\bf x}, \dot{\bf x},t)$
where ${\bf x} \in \mathbb{R}^n$ and $f$ is jointly convex in both ${\bf x}$ and $\dot{\bf x}$ \cite{aecon,Troutman}.
(This idea also generalizes to non-Euclidean spaces with the notions of geodesic convexity.)

Note that the theorem presented in the previous section does not fall neatly into either of these categories since
no curvature or convexity restrictions are placed on $\mathfrak{g}(x)$. Moreover, as will be shown shortly,
higher dimensional globally optimal solutions that build on the results of the previous section can be constructed which neither correspond to a geodesic in a negatively curved space, nor correspond to $f(\cdot)$ being convex. (Recall that the only constraints to do the operations in the proof were that $\mathfrak{g}(x(t)) > 0$ for all $t \in [0,1]$ and $\mathfrak{g}(x)$ needed to be differentiable.)

One way to construct higher dimensional variational problems with globally optimal solutions is when
$$ \dot{\bf x}^T G({\bf x})\, \dot{\bf x} = \sum_{i=1}^{n} \mathfrak{g}_i(x_i) \,\dot{x}_i^2 . $$
Then $n$ decoupled one-dimensional problems of the sort in the previous section result. Hence, if
a problem exists that can be decoupled into the above form by a change of coordinates,
the variational problem
will have a globally minimal solution
even though $G({\bf x})$ may not correspond to a space of negative curvature, nor would $f({\bf x}, \dot{\bf x},t)$ necessarily be convex.

The following theorem addresses a class of multi-dimensional globally optimal variational problems which will be useful in the context
of joint optimization over an externally acting nuisance group and the internally acting TRG. In this theorem, the notation
$\|{\bf s}\|_W = ({\bf s}^T W {\bf s})^{\frac{1}{2}}$ (the weighted Euclidean vector norm with
symmetric positive definite $m\times m$ matrix $W$) is used.
\\ \\
\noindent
{\bf THEOREM 2}:
{\it Suppose that the Euler Lagrange equations provide a global minimum to the problem in (\ref{Jdefapp31245}) with ${\bf x} \in \mathbb{R}^n$ with specified boundary conditions ${\bf x}(0)$ and ${\bf x}(1)$.
Then if ${\bf s} \in \mathbb{R}^m$, the new variational problem in the variable ${\bf q} = [{\bf x}^T, {\bf s}^T]^T \in \mathbb{R}^{n+m}$
will have a globally optimal solution with specified boundary conditions ${\bf q}(0)$ and ${\bf q}(1)$
when $f$ is replaced with $\phi$ as
\begin{equation}
\phi({\bf q}, \dot{\bf q},t) \,\doteq\,
f({\bf x}, \dot{\bf x},t) + c({\bf x},\dot{\bf x},\dot{\bf s},t)
\label{glob2}
\end{equation}
where
$$ c({\bf x},\dot{\bf x},\dot{\bf s},t) \,=\, \frac{1}{2} \, \|\dot{\bf s} - A({\bf x}) \dot{\bf x} \|_W^2, $$
$A({\bf x}) = [A_{ij}{\bf x})]$ is any differentiable $m\times n$ matrix satisfying
\begin{equation}
\frac{\partial A_{ji}}{\partial x_k} = \frac{\partial A_{jk}}{\partial x_i}\,,
\label{acond1}
\end{equation}
and $W = W(t)$ is any differentiable positive definite $m\times m$ matrix function of time.}
\begin{proof}
The Euler-Lagrange equations for this problem,
\begin{equation}
\frac{\partial \phi}{\partial {\bf q}} -
\frac{d}{dt} \left(\frac{\partial \phi}{\partial \dot{\bf q}} \right) = {\bf 0}
\label{ELphi}
\end{equation}
reduce to two sets of equations, one associated with the variable ${\bf x}$ of the form
\begin{eqnarray}
&& \frac{\partial f}{\partial {\bf x}} -
\frac{d}{dt} \left(\frac{\partial f}{\partial \dot{\bf x}} \right)
 - \frac{d}{dt}\left[A^T(\qq) W(\dot{\bf s} - A(\qq)\dot{\qq})\right] \nonumber \\
&&+ \frac{\partial}{\partial \qq}\left[\dot{\qq}^T A^T(\qq)\right]
W(\dot{\bf s} - A(\qq)\dot{\qq}) = {\bf 0},
\label{first3kn2nd}
\end{eqnarray}
and one in ${\bf s}$ of the form
\begin{equation}
\frac{\partial c}{\partial {\bf s}} -
\frac{d}{dt} \left(\frac{\partial c}{\partial \dot{\bf s}} \right) =
\frac{d}{dt} \left\{W(\dot{\bf s} - A({\bf x}) \dot{\bf x})\right\}
= {\bf 0}.
\label{first4kn2nd0}
\end{equation}

Integrating (\ref{first4kn2nd0}) with respect to time gives
\begin{equation}
W(\dot{\bf s} - A({\bf x}) \dot{\bf x}) = {\bf a}
\label{first4kn2nd}
\end{equation}
where ${\bf a}$ is an arbitrary constant vector in $\mathbb{R}^m$.

Substituting (\ref{first4kn2nd}) back into (\ref{first3kn2nd}) and using
the chain rule together with (\ref{acond1}) gives
$$ \frac{d}{dt}\left(A^T(\qq)\right) 
= \frac{\partial}{\partial \qq}\left[\dot{\qq}^T A^T(\qq)\right], $$
which means that (\ref{first3kn2nd}) reduces to (\ref{Jdefapp31245}), and so the
optimal solution for the ``$\qq$-part'' of the problem again will be $\qq^*(t)$ of the original variational
problem in ${\bf x}$. Then, with this $\qq^*(t)$ computed, the solution to (\ref{first4kn2nd0}), or equivalently
(\ref{first4kn2nd}), will be
$$ {\bf s}^*(t) = {\bf b} + \int_{0}^{t} \left\{[W(t')]^{-1} {\bf a} +  A({\bf x}^*(t'))\, \dot{\bf x}^*(t')\right\} \, dt'  $$
where ${\bf a}$ and ${\bf b}$ are determined by fixing ${\bf s}(0)$ and ${\bf s}(1)$.
The cost associated with $f({\bf x}^*, \dot{\bf x}^*,t)$ is as low as it can be since ${\bf x}^*(t)$ is by definition
the global minimizer of the original variational calculus problem. The cost
$$ c({\bf x}^*, \dot{\bf x}^*, \dot{\bf s}^*,t) \,=\, \frac{1}{2} {\bf a}^T W^{-1} {\bf a} $$
is as low as it can be while satisfying initial conditions. This can be observed by adding any perturbation
to ${\bf s}^*(t)$ that preserves the boundary conditions -- the result is an increase in the cost $c$.
Hence ${\bf q}^*(t)$ is the globally optimal solution defined by $({\bf x}^*(t),{\bf s}^*(t))$.
\end{proof}

In the case when
$$ f({\bf x}, \dot{\bf x},t) = \frac{1}{2} \dot{\qq}^T {\cal G}(\qq) \dot{\qq}, $$
then the integrand in the functional for this composite problem can be written as
\begin{equation}
\phi(\qq,\dot{\qq},\dot{\bf s},t) = \frac{1}{2}
\left[\begin{array}{c}
\dot{\qq} \\ \\
\dot{\bf s} \end{array}\right]^T
G(\qq,t)
\left[\begin{array}{c}
\dot{\qq} \\ \\
\dot{\bf s} \end{array}\right]
\label{boot33dd}
\end{equation}
where
$$ G(\qq,t) =
\left(\begin{array}{cc}
{\cal G}(\qq) + A^T(\qq) W(t) A(\qq) & A^T(\qq) W(t) \\ \\
W(t) A(\qq) & W(t) \end{array}\right). $$
This sort of globally optimal variational calculus problem does not generally fall into the negative curvature scenario (even when restricting $W$
to be constant), nor will $f(\cdot)$ be convex in general.

A consequence of this reasoning is that it can be iterated. In other words, now that a globally optimal solution
is obtained to the variational calculus problem with functional $\phi({\bf q}^*(t), \dot{\bf q}^*(t),t)$, an even
higher dimensional problem can be built on this, and so on. This is why the approach is referred to here as bootstrapping.
The next section explains how this theorem can be used in the simultaneous minimization over internal
(temporal) alignment via reparameterization, and external alignment by removal of nuisance group parameters
such as differences in individual perspective, pose, etc.

\section{Searches Over $G\backslash S$} \label{sec2}

Consider a signal that evolves on the intersection of a solid block $\mathbb{B} \subset \mathbb{R}^n$ (i.e., interior of a cube)
and the integer lattice $\mathbb{Z}^n$, resulting in $S=\mathbb{Z}^n \cap \mathbb{B}$, as would be the case for video images when $n=2$. If we think of the $n$-dimensional block of data as being infinitely zero padded, then the block
of data can be viewed as\footnote{Here the values in each matrix entry is taken to be a non-negative
scalar, but for color video it could be viewed as vector valued.}
$$ X_{z_1,z_2,...,z_n} \,=\, h({\bf z}, t) $$
where ${\bf z} = [z_1,z_2,...,z_n]^T \in \mathbb{Z}^n$.
Moreover, we can interpolate pixel values off lattice so that for each value of $t \in [0,1]$ the
function $h: \mathbb{R}^n \times [0,1] \,\rightarrow\, \mathbb{R}_{\geq 0}$ is well defined.

Now suppose that $G$ is a Lie group that reflects nuisance parameters, and suppose that data is collected
from a dynamics scene as
$$ g(t) \cdot X(t) = h([g(t)]^{-1} \cdot {\bf x}, t). $$
(Here the different actions on $X$ and ${\bf x}$ are both denoted as $\cdot$, but there is no ambiguity because
they are clear from the context.) It can be that $g(t)$ is dynamic in scenarios such as a hand-held video camera,
or $g(t)=g_0$ could be a fixed unknown element of $G$. Regardless, it is desirable to quotient out the effects of
$G$.

The static case can be handled by using a metric on the space of images that is left invariant, because then for any pair of two images in the sequence
$$ d(g_0 \cdot X(t_1), g_0 \cdot X(t_1)) = d(X(t_1), X(t_1)). $$
Moreover, in cases where the data evolves directly on $G$, the invariances of $G$ can be
used to quotient out the unknown $g_0 \in G$ and search over a reduced space. An example of this when $G = SE(3)$ is given later in the paper.

In the dynamic scenario, an approach to minimizing the effects of extraneous motion in an image block $X(t)$ is to introduce a variable
$g(t)$ to compensate for unwanted motion by minimizing a cost functional of the form
\begin{equation}
{\it C}_1 \doteq \frac{1}{2}
\int_{0}^{1} \left\{\int_{\mathbb{R}^n} \left|\frac{d}{dt} h([g(t)]^{-1} \cdot {\bf x}, t)\right|^2 d{\bf x}\right\} dt.
\label{nobiggy}
\end{equation}
Using the chain rule while observing that $g = g(t)$ gives
$$ \frac{d}{dt} h(g^{-1} \cdot {\bf x}, t) = $$
$$ [(\nabla h)(g^{-1} \cdot {\bf x}, t)]^T \frac{d}{dt}(g^{-1}) \cdot {\bf x} \,+\, \frac{\partial h}{\partial t} $$
where $\nabla h({\bf x},t) = \partial h/\partial {\bf x}$.
Note that $\frac{d}{dt}(g^{-1}) = - g^{-1} \dot{g} g^{-1} = - \hat{\xi} g^{-1}$ where
$g^{-1} \dot{g} = \hat{\xi} = \sum_k \xi_k E_k$ is the body-fixed velocity associated with $g(t)$ which evolves in the Lie algebra of $G$. $\hat{\xi}$ is expressed in the basis $\{E_i\}$ with $\xi = [\xi_1,...,\xi_N]^T$ where
$N$ is the dimension of $G$. Using this and
making the change of variables ${\bf y} = g^{-1} \cdot {\bf x}$ gives
$$ {\it C}_1 = \frac{1}{2} \int_{0}^{1} \left\{\int_{\mathbb{R}^n} \left|
[(\nabla h)({\bf y}, t)]^T (\hat{\xi}\cdot {\bf y}) - \frac{\partial h}{\partial t} 
\right|^2 d{\bf y}\right\} \Delta(g) dt $$
where the Jacobian determinant $\Delta(g) = |d{\bf x}/d{\bf y}|$ will be equal to unity for groups such as $G=SE(n)$ or $G=SL(n,\mathbb{R})$ acting on $\mathbb{R}^n$ in the usual way, but not in general. For example, for $A \in GL(n,\mathbb{R})$ acting as ${\bf y} = A{\bf x}$, $\Delta(A) = |{\rm det}(A)| \neq 1$.

The structure of the above calculations results in a cost function of the form
\begin{eqnarray}
{\it C}_1 &=& \frac{1}{2} \int_{0}^{1} \left\{\xi^T M(g,t) \xi - 2 \xi^T {\bf b}(g,t) + c(g,t)\right\} \, dt \nonumber \\
&\doteq& \int_{0}^{1} f(g,\xi,t) dt
\label{biggy1}
\end{eqnarray}
Interestingly, in the case when $\Delta(g) = 1$, the quantities $M$, ${\bf b}$, and $c$ become independent of $g$.
Regardless, there are two ways to approach this variational problem. One way would be to introduce coordinates, ${\bf q}$, and
express $g = g({\bf q})$ and $\xi = J({\bf q}) \dot{\bf q}$ and to write
$f(g({\bf q}), J({\bf q}) \dot{\bf q},t) = \phi({\bf q}, \dot{\bf q},t)$ and then to use
(\ref{ELphi}). Alternatively, the lesser-known generalization of the Euler-Lagrange equation
known as the {\it Euler-Poincar\'{e}} equation can directly address variational minimization of (\ref{biggy1}) in a coordinate-free
way by solving
\begin{equation}
\frac{d}{dt}\left(\frac{\partial f}{\partial \xi_i}\right) +
\sum_{j,k=1}^{N} \frac{\partial f}{\partial \xi_k} C_{ij}^{k} \,
\xi_j = \tilde{E}_{i} f
\label{eq:biggy2}
\end{equation}
where $\{E_i\}$ is any basis for the Lie algebra of $G$, the directional derivatives $\tilde{E}_{i} f$ are defined as
$$ (\tilde{E}_{i} f)(g) \doteq \left.\frac{d}{dt}f(g \circ \exp(t E_i))\right|_{t=0}, $$
and $C_{ij}^{k}$ are the structure constants of the Lie algebra such that $[E_i,E_j] = \sum_{k} C_{ij}^{k} E_k$.
For a matrix Lie group such as $SE(3)$, the Lie bracket, $[\cdot,\cdot]$, is simply the matrix commutator.

For a detailed derivation of (\ref{eq:biggy2}) and special cases in which global optimality of solutions to the Euler-Poincar\'{e} equation
can be guaranteed, see \cite{mybook}.

Substituting (\ref{biggy1}) into (\ref{eq:biggy2}) in the case when $\Delta(g) =1$ gives
\begin{eqnarray}
&& \frac{d}{dt}\left(\sum_{l=1}^{N} M_{il}(t) \xi_l - b_i(t) \right) + \nonumber \\
&& \sum_{j,k=1}^{N} \left(\sum_{l=1}^{N} M_{kl}(t) \xi_l - b_k(t)\right) C_{ij}^{k} \, \xi_j = 0.
\label{eq:gfirst}
\end{eqnarray}
If we seek the solution $\xi^*(t)$ that minimizes this with the boundary conditions left free, the result is simply
\begin{equation}
\xi^*(t) = [M(t)]^{-1} {\bf b}(t).
\label{eq:simplymxb}
\end{equation}
But if boundary conditions other than $\xi(0) = [M(0)]^{-1} {\bf b}(0)$ and $\xi(1) = [M(1)]^{-1} {\bf b}(1)$
are required, then (\ref{eq:gfirst}) would need to be solved numerically.

\section{Reducing Searches to $G\backslash S/{\cal T}$} \label{sec3}

Three problems have been considered previously: (1) globally optimal temporal reparameterization; (2) bootstrapping global optimality to higher
dimensional spaces; and (3) variational minimization over Lie groups to ameliorate the effects of nuisance parameters.
This section ties these topics together by addressing simultaneous minimization over $G \times {\cal T}$.
\\ \\
\noindent
{\bf THEOREM 3}:
{\it If $\Delta(g) = 1$, then the globally minimal solution to the variational problem with cost function
$$ {\it C}_2 \doteq \frac{1}{2} \int_{0}^{1} \left\{\int_{\mathbb{R}^n} \left|\frac{d}{dt} h(g^{-1} \cdot {\bf x}, \tau)\right|^2 d{\bf x} \right\} dt $$
in the variables $(g(t),\tau(t)) \in G \times {\cal T}$ with $\tau(0)=0$ and $\tau(1) = 1$ and free boundary conditions on
$g(t)$ and its derivatives, is equivalent to
first globally minimizing the variational problem in (\ref{biggy1})
over $G$, followed by solving the temporal reparameterization problem over ${\cal T}$ in (\ref{optcost}).
}
\begin{proof}
Following the same steps as those which led to (\ref{biggy1}), when $\Delta(g) = 1$, ${\it C}_2$ can be rewritten as
$$ {\it C}_2 = \frac{1}{2} \int_{0}^{1} \left\{\int_{\mathbb{R}^n} \left|
[(\nabla h)({\bf y}, \tau)]^T (\hat{\xi}\cdot {\bf y}) - \frac{\partial h}{\partial \tau} \dot{\tau}
\right|^2 d{\bf y}\right\} dt. $$
The result is of the form
\begin{eqnarray}
{\it C}_2 &=& \frac{1}{2} \int_{0}^{1} \left\{\xi^T M(\tau) \xi - 2 \xi^T {\bf b}(\tau) \dot{\tau}
+ c(\tau) \dot{\tau}^2 \right\} \, dt \nonumber \\
&\doteq& \int_{0}^{1} f(\xi,\tau,\dot{\tau}) \, dt\,,
\label{biggy164}
\end{eqnarray}
where the components of ${\bf b}(\tau)$ are
$$ {b}_i(\tau) = \int_{\mathbb{R}^n}
[(\nabla h)({\bf y}, \tau)]^T (E_i \cdot {\bf y})  \frac{\partial h}{\partial \tau}({\bf y}, \tau) d{\bf y} \,. $$

The resulting variational equations (combination of E-P and E-L) are
\begin{eqnarray}
&& \frac{d}{dt}\left(\sum_{l=1}^{N} M_{il}(\tau) \xi_l - b_i(\tau) \dot{\tau} \right) + \nonumber \\
&& \sum_{j,k=1}^{N} \left(\sum_{l=1}^{N} M_{kl}(\tau) \xi_l - b_k(\tau)\dot{\tau} \right) C_{ij}^{k} \, \xi_j = 0
\label{eq:sametime1}
\end{eqnarray}
and
\begin{equation}
\frac{d}{dt}\{c(\tau) \dot{\tau} - \xi^T {\bf b}(\tau)\} = \frac{1}{2} \xi^T \frac{dM}{d\tau} \xi
- \xi^T \frac{d{\bf b}}{d\tau} \dot{\tau} + \frac{1}{2} \frac{dc}{d\tau} \dot{\tau}^2 .
\label{eq:sametime2}
\end{equation}

\noindent
Denote the solution to the variational problem in (\ref{biggy1})
as $\xi_1^*(t)$, which for given $g(0)$ defines $g_1^*(t)$. Explicitly, $\xi_1^*(t)$ is given in
(\ref{eq:simplymxb}) when the boundary conditions on $\xi(t)$ are free and $\Delta(g) = 1$.
Let $g(t) = g_1^*(\tau(t))$ and use the chain rule. This results in $\xi(t) = \xi_1^*(\tau(t)) \, \dot{\tau}(t)$, which after substituting in (\ref{biggy164}) changes the functional to
\begin{equation}
f(\xi_1^*,\tau,\dot{\tau}) = \mathfrak{g}_2(\tau) \dot{\tau}^2
\label{tau2def}
\end{equation}
where
\begin{eqnarray*}
\mathfrak{g}_2(\tau) \,&=&\, \frac{1}{2} (\xi_1^*)^T M(\tau) \xi_1^* -  (\xi_1^*)^T {\bf b}(\tau)
+ \frac{1}{2} c(\tau) \\
\,&=&\, \frac{1}{2} \{c(\tau) - {\bf b}^T(\tau) [M(\tau)]^{-1} {\bf b}(\tau)\}.
\end{eqnarray*}
The notation $\mathfrak{g}_2(\tau)$ denotes that the computation of this $\mathfrak{g}(\tau)$ follows the computation of $\xi_1^*(t)$.

Variationally minimizing the functional with (\ref{tau2def}) over
over $\tau(t)$, which is of the form in Theorem 1, then gives $\tau_2^*(t)$.
It is easy to see that $\xi = \xi_1^*(\tau) \dot{\tau} $ solves (\ref{eq:sametime1}) for the same reason that
$\xi_1^*(t)$ solves (\ref{biggy1}), independent of the behavior of $\tau(t)$. Moreover, substituting
$$ (g(t), \tau(t)) = (g_1^*(\tau_2^*(t)), \tau_2^*(t)), $$
and hence
\begin{equation}
\xi(t) = \xi_1^*(\tau_2^*(t)) \, \dot{\tau}_2^*(t),
\label{eq:xi1*def}
\end{equation}
into (\ref{eq:sametime2}) reduces to exactly the same thing as (\ref{origode}) with $x(t) = \tau_2^*(t)$ and
$\mathfrak{g}(x) = \mathfrak{g}_2(\tau)$ since
$$ c(\tau) \dot{\tau} - (\xi_1^*)^T {\bf b}(\tau) = 2\,\mathfrak{g}_2(\tau) $$
and
$$ \frac{1}{2} (\xi_1^*)^T \frac{dM}{d\tau} \xi_1^*
- (\xi_1^*)^T \frac{d{\bf b}}{d\tau} + \frac{1}{2} \frac{dc}{d\tau} = \frac{d\mathfrak{g}_2}{d\tau}. $$

Therefore solving (\ref{eq:sametime1}) after first computing $\xi_1^*(\tau) \dot{\tau}$
in for $\xi(t)$ reduces (\ref{eq:sametime1}) to the variational problem of minimizing the
functional with integrand (\ref{tau2def}), the form of which is known from Theorem 1 to produce a global minimum.
\end{proof}


Note that the integrand in (\ref{biggy164}) can be rewritten as
$$ \xi^T\left(M(\tau) - \frac{{\bf b}(\tau) {\bf b}^T(\tau)}{c(\tau)}\right) \xi + c(\tau) \left(\dot{\tau} - \frac{\xi^T {\bf b}(\tau)}{c(\tau)}\right)^2, $$
which becomes the bootstrapped cost in Theorem 2 in the trivial case when $M$, ${\bf b}$, and $c$
are all independent of $\tau$.

\section{An Action Recognition Example}

Consider a random person who is asked to come into a room and stand at an arbitrary position and orientation and is asked to remain
stationary other than moving his/her arm.
Assume that the room is retrofitted with video and/or RGBD cameras so that this imaging system unambiguously recovers a trajectory of
the person's shoulder, elbow, and hand, so that a trajectory: $X_1(t) = (S(t), E(t), H(t)) \in SE(3)\times SE(3)\times SE(3)$ is observed.
The goal is to determine whether or not this trajectory might correspond to a known behavior. Suppose that in a database, trajectories for the acts
of waving, throwing, scratching one's head, giving a thumbs up, raising a hand to ask a question, and rubbing ones' eyes are stored.
Suppose that the trajectories for each of these behaviors has already been observed from several recordings of previous people also at
random positions in the room, and and has already been stored with temporal variations quotiented out, as well as normalizing for scaling effects due to
the different sizes of the people. Let $X_2(t)$ be any of these prior annotated behaviors which have been stored using an optimally
reparameterized timescale.

The methodology presented earlier in this paper then provides a way to rapidly compare $X_1(t)$ with each candidate $X_2(t)$ by allowing
the freedom to temporally reparameterize $X_1(t)$ as $X_1(\tau(t))$. The problem of quotienting out the effects of nuisance groups in this scenario is particularly amenable to efficient solution because $S = [SE(3)]^3$ is a Lie group, which is acted on by $g(t) \in G=SE(3)$ as
$$ g(t) \cdot (S(t), E(t),H(t)) = (g(t)\,S(t), g(t)\,E(t), g(t)\,H(t)). $$
Since both the current person as well as those who contributed to the database stand in the room at random positions and orientations,
one might think that the effect of this unknown pose (which corresponds to $g$ in the above equation and in
(\ref{Ddef})) would need to be accounted for. It might seem like this would require a lot of effort to sample a lot of values to
do the minimization, but there are better ways. For example, new variables $(S(t)^{-1} E(t),$  $S(t)^{-1}H(t),$  $E(t)^{-1}H(t))$ could be recorded
that automatically quotient out the effect of $G=SE(3)$ in (\ref{Ddef})). Similarly, using the trajectory as an object
rather than each instance in the trajectory gives
$$ X_1(t)^{-1} X_1(t + \Delta t) = $$
$$ (S(t)^{-1} S(t + \Delta t), E(t)^{-1} E(t + \Delta t), H(t)^{-1} H(t + \Delta t)). $$
which is invariant to $X$, and which can be optimally reparameterized using the methodology resulting from Theorem 1 to compare
directly with previously observed behaviors.

Another way to minimize the effects of the nuisance group is by using the fact that $SE(3)$ has two screw invariants.
One corresponds to the angle $\theta$ in the expression $R = \exp(\theta \hat{\bf n})$ where
$\hat{\bf n}$ is a $3\times 3$ skew-symmetric matrix with the property that $\hat{\bf n} {\bf x} = {\bf n} \times {\bf x}$ where ${\bf n}$ is the
axis of rotation and ${\bf x} \in \mathbb{R}^3$ is arbitrary. The second, $d = {\bf n}^T {\bf t}$ where ${\bf t}$ is the translation vector in
(\ref{se3def}), corresponds to the distance travelled along the screw axis. The meaning of these invariants is that they are constant with
respect to conjugation in the sense that for any $A,B \in SE(3)$
$$ \theta(A B A^{-1}) = \theta(B) \,\,\, {\rm and} \,\,\, d(A B A^{-1}) = d(B). $$
Consequently, even though $(E(t) S(t)^{-1},$ $H(t) S(t)^{-1},$ $H(t) E(t)^{-1})$ and
$X_1(t + \Delta t)^{-1} X_1(t)$ (with inverses on the right), the screw parameters provide two scalar signals that can be be used
to match with scalar signals in a database. And so reparameterization can take place either at the level of a trajectory in the group $[SE(3)]^3$
with $g(t) \in SE(3)$ quotiented out, or the reparameterization can take place at the level of signals in the space of invariants.

In contrast, if only one rigid-body can be tracked instead of multiple features such as shoulder, elbow and hand, then
there is no way to remove the effects of $g(t)$ by considering relative motions, and then the heavier computation involved in
Theorem 3 becomes the appropriate tool.


\section{Conclusions}

This paper presented a framework for globally optimal reparameterization of the temporal dependence of signals (or trajectories).
This serves as a way to optimally align two signals with $O(N)$ computations when each signal consisting of $N$ values (e.g,
amplitudes, images, feature vectors, etc.). This linear
performance is achieved because the reparameterization depends on the local rate of change of the signal in each individual sequence
and an integral over each individual sequence, followed by pointwise matching at the end of the process when both signals
are renormalized to a universal time scale. The methodology builds on a class of problems in variational calculus with
globally minimal solutions, which appears not to be known in the wider literature.
The resulting method is in contrast to how matching is currently done, which often uses dynamic time warping.
The internal (temporal) reparameterization approach introduced here automatically provides a way to directly compare scalar signals.
Moreover, it is shown that for the case of multidimensional signals, the effects of nuisance groups (such as those due to changes in pose, perspective,
etc.) which act ``externally'' (i.e., at each instant of time) can be eliminated either a priori by using invariants, or simultaneously with the
the temporal reparameterization.
\\ \\
\noindent
{\bf Acknowledgements}
The author thanks Ms. Yixin Gao, Mr. Qianli Ma, Mr. Sipu Ruan, and Profs. Raman Arora and Mauro Maggioni for useful discussions leading to this work.
This work was performed under Office of Naval Research Award N00014-17-1-2142 (GRANT12203025) under code 311.

{\small

\end{document}